\documentclass[a4paper,leqno,twoside,12pt]{article}
\usepackage{amsmath,amssymb,amsthm}
\usepackage{url}

\swapnumbers
\newtheorem{theorem}{Theorem}
\newtheorem{corollary}[theorem]{Corollary}
\newtheorem{lemma}[theorem]{Lemma}
\newtheorem{proposition}[theorem]{Proposition}
\theoremstyle{definition}
\newtheorem{remark}[theorem]{Remark}
\newtheorem{definition}[theorem]{Definition}

\DeclareFontFamily{U}{UWCyr}{}
\DeclareFontShape{U}{UWCyr}{m}{n}{%
  <5> <6> <7> <8> <9>
  <10> <10.95> <12> <14.4> <17.28> <20.74> <24.88> wncyr10
  }{}
\DeclareMathAlphabet{\cyrm}{U}{UWCyr}{m}{n}
\DeclareSymbolFont{cyrm}{U}{UWCyr}{m}{n}
\DeclareSymbolFontAlphabet{\cyrm}{cyrm}
\DeclareMathSymbol{\Evo}{\cyrm}{cyrm}{"03}
\DeclareMathOperator{\byd}{\raisebox{-.2ex}{$\overset{\text{\tiny
def}}{=}$}} \DeclareMathOperator{\im}{im}
 \DeclareMathOperator{\id}{id}

\DeclareMathOperator{\Vol}{\mathcal{V}}
\DeclareMathOperator{\Gr}{Gr} \DeclareMathOperator{\Hess}{Hess}
\newcommand{\cprime}{\/{\mathsurround=0pt$'$}}

\newcommand{\bsi}{\boldsymbol{\sigma}}
\newcommand{\bta}{\boldsymbol{\tau}}

\newcommand{\hd}{\bar{d}}
\newcommand{\R}{\mathbb{R}}

\newcommand{\BL}{\mathbb{L}}
\newcommand{\BT}{\mathbb{T}}
\newcommand{\BM}{\mathbb{M}}
\newcommand{\cC}{\mathcal{C}}
\newcommand{\cE}{\mathcal{E}}
\newcommand{\cF}{\mathcal{F}}
\newcommand{\cH}{\mathcal{H}}
\newcommand{\cL}{\mathcal{L}}

\newcommand{\cA}{\mathcal{A}}

\newcommand{\Diff}{\mathrm{Dif{}f}}
\DeclareMathOperator{\alt}{alt} 
\DeclareMathOperator{\Hom}{Hom} 
\newcommand{\hL}{\Bar{\Lambda}}
\newcommand{\CDiffalt}[2]{\cC\Diff^{\alt}_{(#1)\,#2}}
\newcommand{\CDiff}{\mathcal{C}\mathrm{Diff}}

\newcommand{\odx}{\overline{dx}{}}
\newcommand{\abs}[1]{\lvert#1\rvert}

\newcommand*{\pd}[2]{\mathchoice{\frac{\partial #1}{\partial #2}}
  {\partial #1/\partial #2}{\partial #1/\partial #2} {\partial%
  #1/\partial #2}}
\newcommand{\eval}[2][\right]{\relax
  \ifx#1\right\relax \left.\fi#2#1\rvert}
\let\abs=\envert
\providecommand{\href}[2]{#2} 
\newcommand*{\eprint}[2][]{\href{http://arXiv.org/abs/#2}%
{\begingroup \Url{arXiv:#2}}}
  \title{Geometric aspects of higher order variational principles on
  submanifolds}
 \author{Gianni Manno, Raffaele Vitolo\\
{\footnotesize Dept.\ of Mathematics ``E. De Giorgi'',
University of Lecce}\\
{\footnotesize  via per Arnesano, 73100 Lecce, Italy}\\
{\footnotesize  \texttt{gianni.manno@unile.it},
  \texttt{raffaele.vitolo@unile.it}}
}

\begin{document}

\maketitle

\begin{abstract}
  The geometry of jets of submanifolds is studied, with special interest in
  the relationship with the calculus of variations.  A new intrinsic geometric
  formulation of the variational problem on jets of submanifolds is given.
  Working examples are provided.


  \smallskip

\end{abstract}

\section*{Introduction}

Jets of submanifolds (also known as \emph{manifolds of contact elements}) are
a natural framework for a geometric study of differential equations and the
calculus of variations \cite{AVL91,Ded77,Ded78,Ehr52,KLV86,MoVi94,Vin84}.  The
space of $r$-th order jets of submanifolds $J^r(E,n)$ is introduced through
the notion of contact of order $r$ between $n$-dimensional submanifolds of a
given manifold $E$. These spaces generalize jets of local sections in the
sense that submanifolds which are not transversal to a fibration are also
considered.

In this paper we devote ourselves to the calculus of variations on
$J^r(E,n)$. This subject has been started in a modern framework in
the pioneering papers \cite{Ded78,Vin78}, where the $\cC$-spectral
sequence was introduced (see section~\ref{sec_var_seq}). The main
problem with respect to jets of fiberings is the absence of a
distinguished space of independent variables. This complicates the
computations of the terms of the $\cC$-spectral sequence.

Let us make the above problem more clear with an example. An $r$-th order
Lagrangian on a bundle $\pi\colon E\to M$ is a section $\lambda\colon
J^r\pi\to\wedge^nT^*M$, where $J^r\pi$ is the $r$-th order jet of $\pi$ and
$n=\dim M$. This section can be regarded as an equivalence class
$\lambda=[\alpha]\in\Lambda^n_r/\cC^1\Lambda^n_r$, where $\Lambda^n_r$ is the
space of $n$-forms on $J^r\pi$ and $\cC^1\Lambda^n_r$ is the subspace of
$1$-contact $n$ forms, \emph{i.e.} $n$-forms vanishing on (the $r$-order
prolongation of) sections of $\pi$. A further property of such forms is that
they yield no contribution to action-like functionals. The space
$\Lambda^n_r/\cC^1\Lambda^n_r$ is an element of the first term of the
$\cC$-spectral sequence. The problem is then: how to represent an object in
$\Lambda^n_r/\cC^1\Lambda^n_r$ in view of the absence of a distinguished space
of independent variables?

This problem was first considered in \cite{Ded77}, where it was proposed to use
a sheaf of local Lagrangians whose difference on the intersection of
neighborhoods was a contact form.

In this paper we discuss the use of the bundle of ``truncated''
total derivatives $H^{r+1,r}\to J^{r+1}(E,n)$ as a natural
analogue of the bundle $TM$ in the non-fibered case. We show that
this bundle can be used to represent forms in the first term of
the $\cC$-spectral sequence. Note that the bundle $H^{1,0}$ was
introduced in \cite{Kol73} with the purpose of studying geometric
objects, \emph{i.e.} tensor fields on submanifolds of a given
manifold which depend on derivatives of the immersion, in
homogeneous manifolds. Later on, the bundle $H^{1,0}$ appeared in
\cite{MoVi94} with the purpose of studying higher order analogues
of connections on jets of submanifolds.

We use the bundle $H^{r+1,r}$ to represent Lagrangians, Euler--Lagrange and
Helmholtz morphism so that we provide a geometric (\emph{i.e.}, invariant,
co\-or\-di\-nate\--\-free) formulation of variational problems on jets of
submanifolds. Our formulation reduces to well-known formulations in the case
of jets of fiberings (see for example~\cite{Sau89}), and is a radical
improvement of the old formulation by Dedecker~\cite{Ded77} (see
section~\ref{sec:var_princ} for details).

The theory is illustrated by two examples: the minimal
submanifolds equation and the equation of relativistic particle
motion. In the first example we study the calculus of variations
on geometric objects arising in Riemannian geometry.  An
alternative geometric approach can be found in \cite{BGG02} in the
framework of exterior differential systems.  In the second example
we consider the relativistic mechanics introduced in the framework
of jets of submanifolds in \cite{JaMo96} (see also
\cite{MaVi04,Vit00}) and discuss its variational formulation.

%-----------------------------------------------------------------------------
\section{Jet spaces}\label{sec:jets}

Here we recall the main definitions used in the present paper. Our
main sources are \cite{AVL91,Many99,KLV86,VinII}. By $J^{r}(E,n)$
we denote the $r$-jet of $n$-dimensional (immersed) submanifolds
of an $(n+m)$-dimensional manifold $E$.  We have the obvious
projections
\[
\dots\to J^r(E,n)\overset{\pi_{r,r-1}}{\to}
J^{r-1}(E,n)\to\cdots\to J^{2}(E,n)\overset{\pi_{2,1}}{\to}
J^{1}(E,n) \overset{\pi_{1,0}}{\to} E,
\]
whose inverse limit is the infinite order jet $J^\infty(E,n)$. We
denote by $j_rL\colon L\to J^r(E,n),\,p\mapsto [L]^r_p$, the
prolongation of an $n$-dimensional submanifold $L\subset E$. We
set $L^{(r)}\byd j_{r}L(L)$. \emph{R-planes} are the tangent
planes to $r$-th order prolonged submanifolds. For any $[L]^r_p\in
J^{r}(E,n)$, $R$-planes passing through it biunivocally correspond
to $(r+1)$-jets projecting on $[L]^r_p$: namely, $[L]^{r+1}_{p}$
corresponds to the tangent plane
$R_{[L]^{r+1}_{p}}=T_{[L]^r_p}L^{(r)}$. The \emph{Cartan subspace}
$\cC^r_{\theta}\subset T_\theta J^r(E,n)$ is defined as the span
of all $R$-planes at the point $\theta$, so that we have the
\emph{Cartan
  distribution} $\cC^r\colon\theta\mapsto\mathcal{C}^r_{\theta}$ on $J^r(E,n)$.
We have
\begin{equation}\label{eq.splitting.Cartan}
  \cC^r_{\theta}=R_{\overline{\theta}}\oplus \ker
  d_{\theta}\pi_{r,r-1}, \quad
  \overline{\theta}\in\pi^{-1}_{r+1,r}(\theta).
\end{equation}
A diffeomorphism of $J^r(E,n)$ is called a \emph{contact
transformation} if it preserves $\cC^r$. Any such transformation
can be lifted to a jet space of higher order. A vector field on
$J^{r}(E,n)$ whose local flow consists of contact transformations
is called a \emph{contact field}.
\begin{remark}\label{rem:parametric-approach}
  The space $J^r(E,n)$ can be also defined as the quotient of the higher order
  tangent bundle of regular $(n,r)$-velocities $\operatorname{reg} T^r_nE$ (it
  is the set of $r$-jets of local immersions $\R^n\to E$ at $0$, see
  \cite{KMS93}) with respect to the group $G^r_n$ of $r$-jets of
  `reparametrizations' \cite{Ehr52,Gri97,KMS93,Kru01}.
\end{remark}
{\textbf{Notation:}}
In what follows all manifolds and maps are smooth. Consider a
manifold $E$, $\dim E=n+m$. Greek indexes $\lambda$, $\mu$ run
from $1$ to $n$ and Latin indexes $i$, $j$ run from $1$ to $m$.
Multiindexes will be denoted by the further Greek letters $\bsi$,
$\bta$, where $\bsi=(\sigma_1,\ldots,\sigma_k)$ with
$1\leq\sigma_j\leq n$ and $|\bsi|\byd k\leq r$ (and analogously
for $\bta$).  A \emph{divided chart} on $E$ is a chart of the form
$(x^\lambda, u^i)$, where $1\leq \lambda\leq n$ and $1\leq i\leq
m$.  Unless otherwise specified, a divided chart induces a chart
on $J^r(E,n)$, that we shall denote by $(x^\lambda,u^i_{\bsi})$.
Einstein convention will be used.
%-----------------------------------------------------------------------------%
\subsection{Contact sequences}
\label{ssec:contactseq}

$R$-planes allow us to construct a short exact sequence of vector
bundles over jets. This construction already appeared
in~\cite{Kol73} for the purpose of studying geometric objects on
submanifolds and in \cite{MoVi94} for the study of generalized
connections on jets of submanifolds.

For $r\geq 0$, consider the following bundles over $J^{r+1}(E,n)$: the
pull-back bundle
\begin{equation}
  T^{r+1,r}\byd J^{r+1}(E,n)\underset{J^{r}(E,n)}{\times}T
  J^{r}(E,n),
\end{equation}
the subbundle $H^{r+1,r}$ of $T^{r+1,r}$ defined by
\begin{equation}
  H^{r+1,r}\byd\left\{ \left( [L]_p^{r+1} ,\upsilon \right) \in
    T^{r+1,r} \mid \upsilon \in R_{[L]_p^{r+1}}\right\},
\end{equation}
and the quotient bundle $V^{r+1,r}\byd T^{r+1,r}/H^{r+1,r}$. The
bundles $H^{r+1,r}$ and $V^{r+1,r}$ are strictly related with the
horizontal and vertical bundles in the case of jets of fiberings
(see remark~\ref{re:fibr}).
\begin{definition}
  We call $H^{r+1,r}$ and $V^{r+1,r}$, respectively, the
  \emph{pseudo-horizontal} and the \emph{pseudo-vertical} bundle of $J^r(E,n)$.
\end{definition}
The pseudo-horizontal bundle has some additional features. The following
isomorphism over $\id_{J^{r+1}(E,n)}$ holds:
\begin{equation}\label{eq:1}
  H^{r+1,r} \to J^{r+1}(E,n)\times _{J^{1}(E,n)}H^{1,0}, \quad
  ([L]_{p}^{r+1},\upsilon)\mapsto \left([L] _{p}^{r+1},d\pi_{r,0}(\upsilon)\right).
\end{equation}
The restriction of $H^{r+1,r}$ to an $n$-dimensional submanifold
$L\subset E$ is isomorphic to $TL$, as it is easily seen. By
definition we have the following \emph{contact exact sequence}
\begin{equation}\label{eq:contact}
  0 \to H^{r+1,r} \overset{D}{\hookrightarrow}
  T^{r+1,r} \overset{\omega}{\to} V^{r+1,r} \to  0,
\end{equation}
where $D$ and $\omega$ are, respectively, the natural inclusion and quotient
projection.
\begin{remark}\label{re:fibr}
  If $E$ has a fiber structure $\pi :E\to M$, then $H^{r+1,r}$ is isomorphic
  to $J^{r+1}\pi\times_M TM$ (the horizontal bundle), and $V^{r+1,r}$ to
  $J^{r+1}\pi\times_{J^r\pi}\ker d\pi_r$, where
  $\pi_r=\pi\circ\pi_{1,0}\circ\dots\circ\pi_{r,r-1}$ (the vertical bundle),
  so that the above contact sequence splits.
\end{remark}
Let us evaluate the coordinate expressions of $D$ and $\omega$. A local basis
of the space of sections of the bundle $H^{r+1,r}$ is
\begin{equation*}
  D_\lambda=\frac{\partial }{\partial x^{\lambda }}+
  u_{\bsi,\lambda}^{j}\frac{\partial }{\partial u_{\bsi}^{j}},\quad
  |\bsi|\leq r
\end{equation*}
where the index $\bsi,\lambda$ stands for
$(\sigma_1,\dots,\sigma_s,\lambda)$ with $s\leq r$.  A local basis
of the space of sections $(H^{r+1,r})^*$ dual to $\{D_\lambda\}$
is given by the restriction of the $1$-forms $dx^\lambda$ to
$H^{r+1,r}$, and is denoted by $\overline{dx}{}^\lambda$.  The
local expression of $D$ turns out to be
$D=\overline{dx}{}^\lambda\otimes D_\lambda$.  A local basis of
the space of sections of the bundle $V^{r+1,r}$ is
\[
B^{\bsi}_{j}\byd\left[\frac{\partial}{\partial u^j_{\bsi}} \right]\,\,, \qquad
\abs{\bsi}\leq r.
\]
The local expression of $\omega$ turns out to be $\omega=\omega^j_{\bsi}\otimes
B^{\bsi}_{j}$, where $\omega^j_{\bsi}=
du_{\bsi}^{j}-u_{\bsi,\lambda}^{j}dx^{\lambda}$.

Let $VJ^1(E,n)\byd\ker d\pi_{1,0}$. We have the following
\begin{lemma}\label{le:vert_iso}
  $VJ^1(E,n)\simeq (H^{1,0})^*\otimes_{J^{1}(E,n)}V^{1,0}$.
\end{lemma}
%\begin{proof}
%  Any point of $J^1(E,n)$ can be seen as the inclusion of an $n$-dimensional
%  subspace of $TE$ into $TE$ itself through $D$, hence as a linear map
%  $\overline{dx}{}^\lambda\otimes (\pd{}{x^{\lambda}}+u^i_\lambda \pd{}{u^i})$.
%  A curve tangent to the fiber of $\pi_{1,0}$ at such a point has the tangent
%  vector $\overline{dx}{}^\lambda\otimes \dot u^i_\lambda \pd{}{u^i}$; this
%  proves the above isomorphism.
%\end{proof}
%
The following theorem has been proved in \cite{MMR00}, and in
\cite{MaThesis} with the help of previous lemma.
\begin{theorem}
  For $r\geq 1$, the bundles $\pi_{r+1,r}$ are affine bundles associated with
  the vector bundles $ (\odot^{r+1}(H^{1,0})^*) \otimes_{J^r(E,n)}V^{1,0}$,
  where $\odot$ stands for the symmetric tensor product.
\end{theorem}
In the case $r=0$, $\pi_{1,0}$ coincides with the Grassmann bundle
of $n$-dimensional subspaces in $TE$: in fact $\Gr (T_pE,n)\simeq
\pi_{1,0}^{-1}(p)$. This contrasts with the case of jets of
fiberings. Now, let $(x^\lambda,u^i)$ and $(y^\mu,v^j)$ be two
coordinate charts on $L\subset E$; let us denote by
$(J_\lambda^\mu,J_i^\mu,J_\lambda^j,J_i^j)$ the Jacobian of the
change of coordinates. Then the fibered coordinate change is given
by the following formula
\begin{equation}\label{eq:changefirst}
  v^j_\mu=A_\mu^\lambda(J^j_\lambda+J^j_i u^i_\lambda).
\end{equation}
where $A_\mu^\lambda\byd (J^\mu_\lambda+J^\mu_i u^i_\lambda)^{-1}$.
%-----------------------------------------------------------------------------%
\subsection{Forms and vector fields on jets}
\label{ssec:forms}

Here we study the spaces of forms on finite order jets, and in
particular \emph{contact} and \emph{horizontal} forms, in view of
a geometrical formulation of variational principles.
%Contact forms vanish when restricted to any
%prolonged submanifold; horizontal forms are forms which do not contain contact
%factors.
%We show that horizontal forms have a special polynomial structure
%which is analogue to the case of jets of fiberings \cite{AnDu80,Vit98}. We
%prove that such a structure exists also for horizontal forms on zero-order
%jets.

We denote by $\cF_r$ the algebra $C^{\infty}(J^r(E,n))$.  For $k\geq 0$ we
denote by $\Lambda^{k}_r$ the $\cF_r$-module of $k$-forms on $J^r(E,n)$, and by
$\chi(J^r(E,n))$ the $\cF_r$-module of vector fields on $J^r(E,n)$. We also set
$\Lambda^{*}_r=\bigoplus _{k}\Lambda^{k}_r$.  We introduce the submodule of
$\Lambda_{r}^{k}$ of the \emph{contact forms} of order $r$
\begin{equation*}
  \cC^{1}\Lambda _{r}^{k}\byd\{\,\alpha \in \Lambda _{r}^{k} \mid
  (j_{r}L)^{\ast }\alpha =0\quad\text{for each $n$-dim. submanifold
  }L\subset E\,\}.
\end{equation*}
We set $\cC^1\Lambda^*_r=\bigoplus_k \cC^1\Lambda^k_r$. Moreover, we define
$\cC^p\Lambda^*_r$ as the $p$-th exterior power of $\cC^1\Lambda^*_r$. Of
course, $\cC^p\Lambda^{k}_r=\cC^p\Lambda^*_r\cap\Lambda^k_r$.
We also consider the $\cF_{r+1}$-module $\Lambda^k_{r+1,r}$ of sections of the
bundle $\bigwedge^k(T^{r+1,r})^*$, the $\cF_{r+1}$-module $\cH^k_{r+1,r}$ of
\emph{pseudo-horizontal} $k$-forms, \emph{i.e.}, sections of the bundle
$\bigwedge^k(H^{r+1,r})^*$, and the $\cF_{r+1}$-module of
\emph{pseudo-vertical} fields $\mathcal{V}^{r+1,r}$, \emph{i.e.}, sections of
$V^{r+1,r}$.
\begin{definition}\label{def:horizon}
  Let $q\in\mathbb{N}$. \emph{Horizontalisation} is the map
  \begin{equation*}
    h^{0,q}\colon   \Lambda _{r}^{q}  \to \cH_{r+1,r}^{q}, \quad
    \alpha \mapsto \wedge^q(D^*)\circ (\pi _{r+1,r}^*\alpha).
  \end{equation*}
  Dually, \emph{verticalization} is the map
  $v\colon\chi(J^r(E,n))\to \mathcal{V}^{r+1,r}$, $v(X)\byd \omega\circ
  X\circ\pi_{r+1,r}$.
\end{definition}
If $\alpha\in\Lambda^q_r$, then we have the coordinate expression
\begin{equation}\label{eq:3}
  \alpha = \alpha
  % indexes of \alpha :
  {_{i_1 \dots i_{h} }^{\bsi_1 \dots \bsi_{h}}}
  {_{\lambda_{h+1} \dots \lambda_{q}}}
  % base of forms :
  du^{i_1}_{\bsi_1}\wedge\dots\wedge du^{i_{h}}_{\bsi_{h}}\wedge
  dx^{\lambda_{h+1}} \wedge\dots\wedge dx^{\lambda_{q}},
\end{equation}
where $0 \leq h \leq q$ and $0\leq \abs{\bsi}\leq r$. Hence
\begin{equation}\label{eq:4}
  h^{0,q}(\alpha) = u^{i_1}_{\bsi_{1},\lambda_{1}} \dots
  u^{i_{h}}_{\bsi_{h},\lambda_{h}} \alpha
  % indexes of \alpha :
  {_{i_1 \dots i_{h} }^{\bsi_1 \dots \bsi_{h}}}
  {_{\lambda_{h+1} \dots \lambda_{q}}}
  % base of forms :
  \overline{dx}{}^{\lambda_{1}} \wedge\dots\wedge \overline{dx}{}^{\lambda_{q}}.
\end{equation}
If $X=a^\lambda\pd{}{x^\lambda} +
b^j_{\bsi}\pd{}{u^j_{\bsi}}\in\chi(J^r(E,n))$, then
$v(X)=(b^j_{\bsi}-a^\lambda u^j_{\bsi,\lambda})B^{\bsi}_j$.

Let us introduce the $\cF_r$-module $\hL_{r}^{q}\byd\im h^{0,q}$.
\begin{theorem}\label{th:hyperjac}
  The space $\hL_{r}^{q}$ consists of elements of
  $\cH^q_{r+1,r}$ whose coefficients are fibered polynomials of degree $q$ in
  the highest order variables (\emph{i.e.}, $u^{i_k}_{\bsi_k,\lambda_k}$ with
  $\abs{\bsi_k}=r$).
\end{theorem}
\begin{proof}
  The property of the statement holds in any chart from the coordinate
  expression (\ref{eq:4}). It remains to prove that this is an intrinsic
  property.

  For $r\geq 1$ the property is intrinsic due to the affine structure of
  $\pi_{r+1,r}$.

  For $r=0$, in a similar way to \eqref{eq:changefirst}, we realize that
  $\overline{dy}^\mu = (J^\mu_\lambda+J^\mu_i
  u^i_\lambda)\overline{dx}^\lambda$. Combining this formula with \eqref{eq:4}
  we deduce that the set of sections of the bundle $\hL_{0}^{q}$
  admit a subspace of sections with polynomial coefficients, even if
  $\pi_{1,0}$ \emph{is not} an affine bundle.
\end{proof}
Note that $\hL_{r}^{q}$ does not coincide with the space of
\emph{all} such polynomial forms unless $n=1$. In fact,
\eqref{eq:4} shows that the coefficients of monomials are
skew-symmetric w.r.t. the exchange of pairs ${}^{i_j}_{\bsi_j}$
and ${}^{i_k}_{\bsi_k}$. This property appears in an analogous way
in the case of jets of fiberings, see~\cite{AnDu80}. Since $D$ is
the identity on $TL^{(r)}$, we have the following
\begin{lemma}\label{lem:app}
  Let $\alpha\in\Lambda^q_r$, with $0\leq q\leq n$. We have $
  (j_{r}L)^{\ast}(\alpha)=(j_{r+1}L)^{\ast}(h^{0,q}(\alpha ))$.
\end{lemma}
As an obvious consequence of the previous lemma, we have:
\begin{equation}\label{cC_and_ker_h}
  \cC^1\Lambda^q_{r} = \ker h^{0,q} \quad \text{if} \quad 0\leq q \leq
  n, \qquad \cC^1\Lambda^q_{r} = \Lambda^q_{r} \quad \text{if} \quad
  q > n.
\end{equation}
Moreover, $\alpha\in\cC^1\Lambda^q_{r}$ if and only if
$\pi_{r+1,r}^*(\alpha)\in\im ((\omega)^*\wedge\,\id)$. Hence, if $\alpha \in
\cC^p\Lambda^{p+q}_{r}$, then
\begin{equation}\label{eq:coord_contact}
  \pi^*_{r+1,r}(\alpha)=\omega^{i_1}_{\bsi_1}\wedge\dots
  \wedge\omega^{i_p}_{\bsi_p}\wedge \alpha_{i_1\dots
    i_p}^{\bsi_1\dots\bsi_p},\qquad \alpha_{i_1\dots
    i_p}^{\bsi_1\dots\bsi_p}\in\pi_{r+1,r}^*\left(\Lambda^q_r\right),
\end{equation}
where $\abs{\bsi_l}\leq r$ for $l=1,\dots,p$. Note that, if $q=0$, then
$\abs{\bsi_l}\leq r-1$ for $l=1,\dots,p$. Moreover, derivatives of order $r+1$
appear in the above expression when $\abs{\bsi_l}=r$. It is possible to obtain
an expression containing just $r$-th order derivatives by using contact forms
of the type $d\omega^{i_l}_{\bsi_l}$ with $\abs{\bsi_l}=r-1$ (see
\cite{Kru90}).

For future purposes, we introduce the following \emph{partial
horizontalisation} map
\begin{equation}\label{eq:5}
  h^{p,q}\colon \Lambda ^{p+q}_r \to \Lambda ^{p}_{r+1,r}\otimes
  \hL^{q}_r,\quad \alpha \mapsto
  (\wedge^p\id\otimes\wedge^q {D}^*)\circ (\pi _{r+1,r}^*\alpha).
\end{equation}
The action of $h^{p,q}$ on decomposable forms is
\begin{multline*}
  h^{p,q}(\alpha_{1}\wedge \ldots \wedge \alpha _{p+q})=
  \\
  \frac{1}{p!\,q!}\sum_{\varsigma \in S_{p+q}}\abs{\varsigma} \pi_{r+1,r}^*(\alpha
  _{\varsigma(1)}\wedge \ldots \wedge \alpha _{\varsigma (p)}) \otimes
  h^{0,q}(\alpha_{\varsigma (p+1)}\wedge \ldots \wedge \alpha _{\varsigma
    (p+q)})\notag ,
\end{multline*}
where $S_{p+q}$ is the set of permutations of $p+q$ elements.

\section{$\cC$-spectral sequence and variational sequence}\label{sec_var_seq}

In this section we show that the terms of the $\cC$-spectral sequence can be
explicitly computed through the pseudo-horizontal and pseudo-vertical bundles
and the horizontalization.

Let $r>0$. We have the \emph{bounded} filtration of modules
\begin{equation}\label{filtration}
  \Lambda^k_r \byd \cC^0\Lambda^{k}_r\supset \cC^1\Lambda^{k}_r \supset
  \dots \supset
  \cC^p\Lambda^{k}_r \supset \dots \supset \cC^I\Lambda^{k}_r \supset
  \cC^{I+1}\Lambda^{k}_r = \{0\}
\end{equation}
(see \cite{Kru90} for the value of $I$).  This is a differential
filtration, \emph{i.e.}
$d(\cC^p\Lambda^{k}_r)\subset\cC^p\Lambda^{k+1}_r$, hence it gives
rise to a spectral sequence $(E^{p,q}_k,d^{p,q}_k)$, with
$k,p,q\geq 0$, which we call the \emph{$\cC$-spectral sequence of
(finite) order $r$} on $E$. By construction, it is invariant
w.r.t.\ contact transformations of $J^r(E,n)$.

We recall that $E_{0}^{p,q}\equiv \cC^p\Lambda^{p+q}_{r} \big /
\cC^{p+1}\Lambda^{p+q}_{r}$. We have the following simple result.
\begin{proposition}
  The restriction of $h^{p,q}$ to $\cC^{p}\Lambda_{r}^{p+q}$ yields the
  injective morphism
  \begin{displaymath}
    E_{0}^{p,q}\to \Lambda_{r+1,r}^{p}\otimes \hL_{r}^{q},\,\quad \ [\alpha]
    \mapsto h^{p,q}(\alpha).
  \end{displaymath}
  It follows that $\hd(h^{p,q}(\alpha ))=h^{p,q+1}(d\alpha)$, where $\hd\byd
  d^{p,q}_0$.
\end{proposition}

The finite order $\cC$-spectral sequence is filtered with respect
to the order $r$: pull-back through $\pi_{r+1,r}$ provides an
injective morphism between the $\cC$-spectral sequences of order
$r$ and $r+1$, respectively. It follows that the infinite order
$\cC$-spectral sequence \cite{Many99,Vin78,Vin84} can be obtained
as the direct limit of the finite order $\cC$-spectral sequences.
Until the end of this section we use the infinite order
formulation for the sake of simplicity. Below we denote by $\cF$
the direct limit of $\cF_r$, and $\hL^{n}$ the direct limit of
$\hL^{n}_r$. Note that $\alpha\in \hL^{n}$ has the coordinate
expression $\lambda=\lambda_0\Vol_n$, where $\lambda_0\in\cF_r$
for some $r$ and $\Vol_n\byd n!\,\odx^1\wedge\dots\wedge\odx^n$ is
a local volume form on any submanifold.

For $0\leq q\leq n$ and $p\geq 1$ the terms $E^{0,q}_0$ and $E^{p,n}_1$ may be
arranged into a further complex which is of fundamental importance for the
calculus of variations (see also the end of section \ref{sec:var_princ}): the
variational sequence.
\begin{definition}\label{def:varseq}
  The complex
    \begin{displaymath}
      \dots \overset{\hd}{\to} \hL^{n-1} \overset{\hd}{\to} \hL^{n}
      \overset{\cE}{\to} E_0^{1,n} /\hd(E_{0}^{1,n-1}) \overset{\cH}{\to}
      E_0^{2,n} /\hd(E_{0}^{2,n-1}) \overset{d^{2,n}_1}{\to} \dots,
  \end{displaymath}
  where $\cE$ is the composition of the projection $\hL^{0,n} \to
  E^{0,n}_1=\hL^{0,n}/\hd(\hL^{0,n-1})$ with $d^{0,n}_1$ and $\cH\byd
  d^{1,n}_1$, is said to be the \emph{variational sequence}.  Elements of
  $E_0^{p,n} /\hd(E_{0}^{p,n-1})$ are said to be \emph{variational $p$-forms}.
\end{definition}
It is easy to prove that $\hd\bar{\alpha}=D_\lambda(\alpha^{\bta_1\cdots\bta_l}
_{j_1\cdots j_l\ \lambda_{l+1}\cdots\lambda_q})
u^{j_1}_{\bta_1,\lambda_1}\cdots u^{j_l}_{\bta_l,\lambda_l}
\odx^{\lambda}\wedge\odx{}^{\lambda_1} \wedge\cdots\wedge\odx{}^{\lambda_{q}}$
for $\abs{\bta_h}=r$. This expression is similar to that of the horizontal
differential in jets of fiberings, but the presence of $\odx{}^\lambda$ yields
different transformation properties.

\smallskip

Now we shall prove that any variational form can be represented by a
distinguished object. To this purpose, the most important tool is Green's
formula \cite{Vin77}, which is the geometric analogue of the integration by
parts.

Let $P$, $Q$ be modules of local sections of vector bundles on
$J^\infty(E,n)$. We recall that a
\emph{$\mathcal{C}$-dif\-fer\-en\-tial operator} $\Delta\colon
P\to Q$ is a linear differential operator which admits a
restriction on prolonged submanifolds $L^{(\infty)}$.  In
coordinates, we have $\Delta=(a_{ij}^{\bsi}D_{\bsi})$, where
$a_{ij}^{\bsi}\in \cF$, $D_{\bsi}=D_{\sigma _{1}}\circ \dots \circ
D_{\sigma _{h}}$, and $\abs{\bsi}=h\leq k$.  The space of such
operators is denoted by $\CDiff_k(P,Q)$. We shall also deal with
spaces of antisymmetric $\mathcal{C}$-differential operators of
$l$ arguments in $P$, which we denote by
$\CDiff^{\alt}_{(l)\,k}(P,Q)$.

Now, let $\hat{P}\byd\Hom(P,\hL^n)$; we say $\hat{P}$ is the
\emph{variational dual} of $P$. Note that the variational dual is
not exactly the same as in \cite{Vin84} because we use our
`concrete' space $\hL^n$ instead of just the quotient
$\Lambda^n/\cC^1\Lambda^n$.

Let $\Delta\in \CDiff(P,Q)$. Then, according with Green's formula
\cite{Vin77,Many99} (see \cite[p.\ 31]{KrVe98} for more details) there exists a
unique $\Delta^*\in \CDiff(\hat{Q},\hat{P})$ such that
\begin{equation}\label{eq:adjoint}
  \widehat{q}(\Delta(p))-(\Delta^*(\widehat{q}))(p) =
  \bar{d}\omega_{p,\widehat{q}}(\Delta)
\end{equation}
for all $\widehat{q}\in\widehat{Q}$, $p\in P$, where
$\omega_{p,\widehat{q}}(\Delta)\in\hL^{n-1}$.  In coordinates, if
$\Delta(\varphi)^i=\Delta_{j}^{\bsi i}D_{\bsi}\varphi^j$, then
$\Delta^*(\psi)_j=(-1)^{\abs{\bsi}}D_{\bsi}(\Delta_{j}^{\bsi i}\psi_i)$.

Let us introduce the $\cF$-module $\varkappa$ of local sections of
$J^\infty(E,n)\times_{J^1(E,n)}V^{1,0}$. With an element
$\varphi\in\varkappa$, a non-trivial symmetry $\Evo_\varphi$ of
the Cartan distribution on $J^\infty(E,n)$ (\emph{i.e.} a symmetry
which is not contained in the distribution) is associated
\cite{Many99}. Locally, if $\varphi=\varphi^iB_i$, then
$\Evo_{\varphi}=D_{\bsi}(\varphi^i)B^{\bsi}_i$. Moreover, let
$K_{p}(\varkappa)\subset
\CDiffalt{p-1}{}(\varkappa,\hat{\varkappa})$ be the subspace of
operators which are skew-adjoint in each argument.  With a
reasoning similar to the fibered case \cite[p.\ 190]{Many99} we
obtain the isomorphisms
\begin{align*}
  & \nabla\colon E^{p,q}_0\to\CDiffalt{p}{r}(\varkappa,\bar{\Lambda}^{q}),
  \quad \nabla_{\beta}\big(\varphi_1,\ldots ,\varphi_p\big)=
  \frac{1}{p!}\; \Evo_{\varphi_1}\lrcorner\left(\cdots \lrcorner
    \Evo_{\varphi_p} \lrcorner \beta\cdots \right)
  \\
  & I_p\colon E^{p,n}_0\big / \hd E^{p,n-1}_0 \to K_{p}(\varkappa), \quad
  I_p([\Delta])(\varphi_{1},\ldots,\varphi_{p-1})(\varphi_p)\byd
  \Delta(\varphi_{1},\ldots,\varphi_{p-1},\cdot)^*(1)(\varphi_p),
\end{align*}
where we set $\beta=h^{p,q}(\alpha)$ and $\Delta=\nabla_\beta$ for
$\alpha\in\cC^p\Lambda^{p+q}$. We have the coordinate expressions
\begin{align}\label{eq:2}
  & \Delta(\varphi_1,\ldots,\varphi_p) =
  \Delta_{i_1\dots i_{p-1}\,j}^{\bsi_1\dots\bsi_{p-1}\bta}
  D_{\bsi_1}(\varphi^{i_1}_1)\cdots
  D_{\bsi_{p-1}}(\varphi^{i_{p-1}}_{p-1})
  D_{\bta}(\varphi^j_p)\,\Vol_n,
  \\
  & I_p([\Delta])(\varphi_1,\ldots,\varphi_p) =
  (-1)^{\abs{\bta}}D_{\bta}(\Delta_{i_1\dots i_{p-1}\,j}^{\bsi_1\dots\bsi_{p-1}\bta}
  D_{\bsi_1}(\varphi^{i_1}_1)\cdots
  D_{\bsi_{p-1}}(\varphi^{i_{p-1}}_{p-1}))
  \varphi^j_p\,\Vol_n,
\end{align}
where $0\leq\abs{\bsi_1},\ldots,\abs{\bsi_{p-1}},\abs{\bta}\leq r$. Note that
if $\alpha\in\cC^p\Lambda^{p+q}_r$ then $\Delta_{i_1\dots
  i_{p-1}\,j}^{\bsi_1\dots\bsi_{p-1}\bta}$ is an $n$-th degree polynomial in
the highest order derivatives with the structure described in
theorem~\ref{th:hyperjac}. A coordinate expression of $\cE$ can be
obtained as follows. Let $\lambda\in\hL^n$, with
$\lambda=\lambda_0\Vol_n$; then the form $\alpha=\lambda_0 n!\,
dx^1\wedge\cdots\wedge dx^n$ fulfills $h^{0,n}(\alpha)=\lambda$,
hence, by definition,
$d^{0,n}_1([\lambda])=([[d\alpha]])=[h^{1,n}(d\alpha)]$. So,
\begin{equation}\label{eq:E-L}
  \cE(\lambda)(\varphi_1)=
  (-1)^{\abs{\bta}}D_{\bta}\left(\pd{}{u^i_{\bta}}\lambda_0\right)
  \, \varphi^i_1\Vol_n.
\end{equation}
In an analogous way we obtain the coordinate expression for the case $p>1$.

%---------------------------------------------------------------------------%

\section{Variational principles}
\label{sec:var_princ}

In this section we formulate higher-order variational principles
on submanifolds in the `classical' way, \emph{i.e.}, with an
integral formalism. The main difference with the fibered case is
the absence of a space of independent variables. We get rid of
this difficulty by the pseudo-horizontal bundle and the
horizontalization operator. A similar approach was already
attempted in \cite{Ded77} where the author was forced to use
families of Lagrangians defined on open subsets with the property
that, on intersecting subsets, the action be the same. In what
follows, the pseudo-horizontal bundle allows us to use single
objects as Lagrangians. More precisely, a Dedecker family of
Lagrangians yields just one bundle morphism $\lambda\colon
J^r(E,n)\to\wedge^n(H^{1,0})^*$. The Euler--Lagrange equations
follow straightforwardly, as we will see.

\begin{definition}\label{def:lag-action}
  A form $h^{0,n}(\alpha)=\lambda \in\hL_r^{0,n}$ is said to be an \emph{$r$-th
    order (generalized) Lagrangian}. The \emph{action} of the Lagrangian
  $\lambda$ on an $n$-dimensional oriented submanifold $L\subset E$ with
  compact closure and regular boundary is defined by
  \begin{equation}\label{eq:11}
    \int_{L}(j_{r}L)^{\ast}\alpha.
  \end{equation}
\end{definition}
The word `generalized' comes from the fact that $\lambda$ depends
on $(r+1)$-st derivatives in the way specified in
equation~\eqref{eq:4}.  Moreover, the action is well-defined
because only the horizontal part of a form $\alpha$ contributes to
the action (lemma~\ref{lem:app}). Now, we formulate the
variational problem. Let $L\subset E$ be as in the above
definition.  A vector field $X$ on $E$ vanishing on $\partial L$
is called a \emph{variation field}.  The submanifold $L$ is
\emph{critical} if for each variation field $X$ with flow $\phi_t$
we have
\begin{equation}\label{eq:extreme_action}
  \frac{d}{dt}\Big |_{t=0} \int_{L}(\phi^{(r)}_t\circ
  j_{r}L)^{\ast}\alpha =0
\end{equation}
where we recall that $\phi^{(r)}_t\colon J^{r}(E,n)\longrightarrow
J^{r}(E,n)$ is the $r$-lift of $\phi_{t}$ (see
section~\ref{sec:jets}). In what follows by $X^{(r)}$ we denote
the $r$-lift of a vector field $X$ on $E$.
\begin{lemma}\label{lemma.h.commuta.con.inserz}
  For any $\alpha\in\cC^1\Lambda^{q+1}$, we have that
  %\begin{equation*}
   $ h^{0,q}(X^{(r)}\lrcorner\alpha) =
    v(X^{(r)})\lrcorner(h^{1,q}(\alpha))$.
  %\end{equation*}
\end{lemma}
\begin{proof}
  In fact,
  \begin{align*}
    h^{0,q}\left(X^{(r)}\lrcorner\alpha\right) &= \wedge^q
    D^{*}\left(\pi^*_{r+1,r}\left(X^{(r)}\lrcorner\alpha\right)\right)
    = \wedge^q
    D^{*}\left(X^{(r)}\lrcorner\left(\pi^*_{r+1,r}\alpha\right)\right)
    \\
    & = \wedge^q D^{*}
    \left(X^{(r)}\lrcorner\left(\omega^i_{\bsi}\wedge\beta^{\bsi}_i\right)\right)
    = \wedge^q
    D^{*}\left(\omega^i_{\bsi}(v(X^{(r)}))\beta^{\bsi}_i\right)
    \\
    & = \left(\left(\id\otimes \wedge^q
        D^{*}\right)(\omega^i_{\bsi}\wedge\beta^{\bsi}_i)\right)v(X^{(r)})
    = v(X^{(r)})\lrcorner\left(h^{1,q}\alpha\right).
  \end{align*}

\end{proof}
\begin{lemma}\label{lemma.magical.Cartan.horiz.formula}
  For any $\beta\in E^{1,q}_0$, we have that
  %\begin{equation*}
   $ v(X^{(r)})\lrcorner(\bar{d}\beta) + \bar{d}(v(X^{(r)})\lrcorner\beta) =
   0$.
  %\end{equation*}
\end{lemma}
\begin{proof}
  Let $\alpha\in\cC^1\Lambda^{1+q}_r$ be such that $h^{1,q}(\alpha)=\beta$. We
  have that
  \[
  v(X^{(r)})\lrcorner\left(\bar{d}(h^{1,q}(\alpha))\right) =
  {v(X^{(r)})}\lrcorner\left(h^{1,q+1}(d\alpha)\right)
  \]
  and in view of lemma~\ref{lemma.h.commuta.con.inserz}, the above right side
  term is equal to $h^{0,q+1}\left({X^{(r)}}\lrcorner d\alpha\right)$.
  Furthermore, taking into account that $\cL_{X^{(r)}}\alpha={X^{(r)}}\lrcorner
  d\alpha+d({X^{(r)}}\lrcorner\alpha)$ is a $1$-contact form, we have that
  \[
  h^{0,q+1}\left({X^{(r)}}\lrcorner d\alpha\right) =
  -h^{0,q+1}\left(d({X^{(r)}}\lrcorner\alpha)\right) =
  -\bar{d}\left(h^{0,q}({X^{(r)}}\lrcorner\alpha)\right) =
  -\bar{d}\left({v(X^{(r)})}\lrcorner h^{1,q}(\alpha)\right)
  \]
  where the last equality is attained by applying again
  lemma~\ref{lemma.h.commuta.con.inserz}.
\end{proof}
\begin{theorem}
  Let $\lambda \in\hL_r^n$. Then an (embedded) submanifold $L\subset E$ is
  critical for $\lambda$ if and only if the following Euler--Lagrange equations
  are fulfilled:
  \begin{equation}\label{eq.Eu-Lag}
  \cE(\lambda)\circ j_{2r}L=0.
  \end{equation}
\end{theorem}
\begin{proof}
  We show that~(\ref{eq:extreme_action}) depends on the vertical part $v(X)$ of
  $X$, and provide the Euler--Lagrange equations. We have
  \begin{multline}
    \frac{d}{dt}\Big |_{t=0}  \int_{L}(\phi^{(r)}\circ j_{r}L)^{\ast}\alpha
    =\int_{L}(j_{r}L)^*\cL_{X^{(r)}}\alpha
    =\int_{L}(j_{r}L)^*{X^{(r)}}\lrcorner d\alpha
        \\
    =\int_L(j_{r+1}L)^*h^{0,n}({X^{(r)}}\lrcorner d\alpha)=
    \int_{L}(j_{r+1}L)^*{v(X^{(r)})}\lrcorner h^{1,n}(d\alpha)
    =\int_{L}(j_{2r+1}L)^*{v(X)}\lrcorner\cE(\lambda).\label{eq:15}
  \end{multline}
  Here we used Stokes' theorem, lemma~\ref{lem:app} and
  lemma~\ref{lemma.h.commuta.con.inserz}, taking into
  account that $d\alpha\in\Lambda^{n+1}_r=\cC^1\Lambda^{n+1}_r$, up to the
  final equality. Green's formula yields
  \begin{equation}\label{eq:quasi_Noether} {v(X^{(r)})}\lrcorner
    h^{1,n}(d\alpha)=\cE(\lambda)(v(X))+\hd
    \omega\,, \quad \omega\in \bar{\Lambda}^{n-1},
  \end{equation}
  hence the last equality of~\eqref{eq:15} follows from the identity
  $(j_{r+1}L)^*\hd\omega=d(j_{r+1}L)^*\omega$, and Stokes' theorem.
  By virtue of the fundamental lemma of calculus of variations,
  equation~\eqref{eq:extreme_action} vanishes if and only if
  $(j_{2r+1}L)^*\cE(\lambda)=0$, which is equivalent to the Euler--Lagrange
  equation \eqref{eq.Eu-Lag}.
\end{proof}
Note that, in view of lemma
~\ref{lemma.magical.Cartan.horiz.formula}, we have
   $ h^{1,n}(d\alpha)=\cE(\lambda)+\hd
    \omega$, with $\omega\in E^{1,n-1}_0$.

It is now clear that in the variational sequence
$E_1^{1,n}$ is the space of \emph{Euler--Lagrange type morphism}, and
$\cE$ takes a Lagrangian $\lambda$ into its Euler--Lagrange form
$\cE(\lambda)$. If $\cE(\lambda)=0$ then the Lagrangian is trivial, or
\emph{null}.

We could continue our analysis and show that $E_1^{2,n}$ is the space of
\emph{Helmholtz type morphism}. The operator $\cH$ takes an Euler--Lagrange
type form $\epsilon$ into its Helmholtz form $\cH(\epsilon)$. If
$\cH(\epsilon)=0$, or, equivalently, if $\epsilon$ is \emph{locally
  variational}, then $\epsilon$ comes from a local Lagrangian.  In other words,
the associated differential equation $\epsilon\circ j_rL = 0$ is locally
variational.

\begin{remark}
  By construction, and in view of Lie--B\"acklund theorem \cite{Many99}, the
  Euler--Lagrange operator $\cE$ and the Helmholtz operator $\cH$ are invariant
  w.r.t.\ contact transformations. In particular, if $m > 1$ then they are
  invariant w.r.t.\ point transformations.
\end{remark}
\begin{remark}
  In the parametric approach~\cite{CrSa03a} the variational principle is
  formulated on $\operatorname{reg} T^r_nE$ under the hypothesis that the
  Lagrangian commute with the action of the group of parametri\-za\-tions
  (remark~\ref{rem:parametric-approach}).  This leads to extra computations in
  order to verify at each step the invariance of objects with respect to
  changes of parametrization.
\end{remark}

%----------------------------------------------------------------------------%
% Section
%----------------------------------------------------------------------------%
\section{Examples}\label{sec:examples}

In this section we will present two examples, one from differential geometry
(the equation of minimal submanifolds) and the other from mathematical physics
(the equation of particle motion in general relativity).

\subsection{The minimal submanifold equation}
\label{sec:minsub}

Here we will show that the geometry of submanifolds of a given
Riemannian manifold can be reformulated on jets of submanifolds.
In fact, tensors that are defined on one submanifold and depend on
derivatives of the immersion are converted into objects that are
defined on jets of submanifolds (\emph{geometric objects}, see
\cite{Kol73} for the case when $E$ is an homogeneous manifold).

Let $(E,g)$ be a Riemannian manifold. Let us denote by $\Gamma$
the Levi--Civita connection associated with $g$. The metric $g$
can be lifted as a fiber metric on $T^{1,0}$ by composition with
the projection $\pi_{1,0}$. We indicate the above metric on
$T^{1,0}$ with the same letter $g$. We also have the contravariant
metric $\bar{g}$ on $(T^{1,0})^*$. Let us set $V^{1,0}_g\byd
(H^{1,0})^\perp$; note that the projection $\omega$ is an
isomorphism between $V^{1,0}_g$ and $V^{1,0}$. We have the
splitting $T^{1,0}=H^{1,0}\oplus V^{1,0}_g$. Note that, if
$L\subset E$ is an $n$-dimensional submanifold of $E$, then the
restriction to $L$ of the bundle $V^{1,0}_g$ is just the normal
bundle $N_gL$ to the submanifold $L$ with respect to the metric
$g$, as it can be easily seen.

The metric $g$ restricts to the metric $g^H$ on $H^{1,0}$, with coordinate
expression
\begin{equation}
  g^H = g^H_{\lambda\mu}\odx^\lambda\otimes\odx^\mu
  = (g_{\lambda\mu} + g_{\lambda j}u^j_{\mu} + g_{i\mu}u^i_{\lambda} +
  g_{ij}u^i_{\lambda}u^j_{\mu})\,\odx^\lambda\otimes\odx^\mu.
\end{equation}
The above metric $g^H$ can be characterized as follows.  If $L\subset E$ is an
$n$-dimensional submanifold of $E$, then $(j^1L)^* g^H$ is the pull-back metric
of $g$ on $L$. For this reason we say $g^H$ is the \emph{universal first
  fundamental form} on $n$-dimensional submanifolds of $E$.

The contravariant metric $\bar{g}^H$ will be used; its coordinate expression is
denoted by $\bar{g}^H = (\bar{g}^H)^{\lambda\mu} D_\lambda\otimes D_\mu$.  We
introduce the local basis $N_i \byd \pd{}{u^i} - (g_{\mu i} + g_{ij} u^j_\mu)
(\bar{g}^H)^{\mu\lambda}D_\lambda$ of $V^{1,0}_g$. The metric $g$ restricts to
the metric $g^V$ on $V^{1,0}_g$, with coordinate expression
\begin{displaymath}
  g^V_{ij} = g(N_i,N_j) = g_{ij} - (g_{\lambda i} +
  g_{ik}u^k_\lambda)(g_{\mu j} + g_{jk}u^k_\mu)(\bar{g}^H)^{\lambda\mu}.
\end{displaymath}
Let us introduce the operator
\begin{equation}
  \textbf{II}\colon J^2(E,n) \to
  (H^{1,0})^*\underset{J^{2}(E,n)}{\otimes}
  (H^{1,0})^*\underset{J^{2}(E,n)}{\otimes }V^{1,0},
  \quad
  \textbf{II}(X,Y,\varrho)=\varrho
  \left(\left[\pi_{2,0}^*(\Gamma)\right]_{\widetilde{X}}(Y)\right)
\end{equation}
where $\widetilde{X}$ is a field on $J^2(E,n)$ lying in the Cartan
distribution and projecting on $X$ (see
\eqref{eq.splitting.Cartan}), and
$[\pi_{2,0}^*(\Gamma)]_{\widetilde{X}}$ is the covariant
derivative, w.r.t. $\widetilde{X}$, of the pull-back connection
$\pi_{2,0}^*(\Gamma)$. The previous definition is well posed as
the vertical part of $\widetilde{X}$ gives no contribution.  We
call \textbf{II} the \emph{universal second fundamental
  form} associated with $\Gamma$. This name is justified by the fact that, if
$L\subset E$ is a submanifold, then the pull-back $(j^2L)^* \textbf{II}\,$
coincides with the second fundamental form on $L$.  Moreover, we call the
following map
\[
\textbf{H}\colon J^2(E,n)\to V^{1,0}_g, \quad \textbf{H}\byd
\frac{1}{n}\, \textbf{II}\circ\bar{g}^H
\]
the \emph{universal mean curvature normal} (or \emph{vector}) on
$n$-dimensional submanifolds of $E$.  It is easy to realize that both
$\textbf{II}^{-1}(0)$ and $\textbf{H}^{-1}(0)$ are regular submanifolds of
$J^2(E,n)$. They are the \emph{totally geodesic submanifold equation} (see
\cite{MaAdvances} for another geometric characterization of this equation) and
the \emph{minimal submanifold equation}.  The coordinate expression of the
latter is
\begin{multline} \label{eq:130}
  \bar{g}^H{^{\lambda\xi}}\left(u^k_{\lambda\xi}+\Gamma_\lambda{}^k{}_\xi +
    \Gamma_\lambda{}^k{}_i \, u^i_\xi + \Gamma_j{}^k{}_\xi \, u^j_\lambda +
    \Gamma_j{}^k{}_i \, u^j_\lambda u^i_\xi \right.
  \\
  \left. -u^k_\beta\left(\Gamma_\lambda{}^\beta{}_\xi +
      \Gamma_\lambda{}^\beta{}_i \, u^i_\xi + \Gamma_j{}^\beta{}_\xi \,
      u^j_\lambda + \Gamma_j{}^\beta{}_i \, u^j_\lambda
      u^i_\xi\right)\right)=0.
\end{multline}
As an example, we write down the equation of minimal surfaces in the Euclidean
space $\R^3$. We suppose that $(x,y,u)$ is a Cartesian coordinate chart. Then
we have
\[
g^H=
\begin{pmatrix}
  1+u_x^2 & u_xu_y\\ u_xu_y & 1+u_y^2
\end{pmatrix},\quad \bar{g}^H=\frac{1}{1+u_x^2+u_y^2}
\begin{pmatrix}
  1+u_y^2 & -u_xu_y\\ -u_xu_y & 1+u_x^2
\end{pmatrix}.
\]
Hence the equation of totally geodesic submanifold is locally
represented by the system $u_{xx}=0$, $u_{xy}=0$, $u_{yy}=0$, and
its contraction with $\bar{g}^H$ yields usual minimal surface
equation $(1+u_y^2)u_{xx}-2u_xu_yu_{xy}+(1+u_x^2)u_{yy}=0$.  This
coordinate representation does not cover the whole equation as a
submanifold of $J^2(\R^3,2)$. Two more divisions of the Cartesian
chart $(x,y,u)$ are needed, namely those for which $x$ and $y$
play the role of dependent variables, respectively.

The most important functional on submanifolds of a Riemannian
manifolds is the area element
  \begin{displaymath}
    \cA \colon J^1(E,n) \to \bigwedge^n (H^{1,0})^*,\quad
    \cA = \sqrt{\abs{g^H}} \, \odx^1\wedge\cdots\wedge \odx^n,
  \end{displaymath}
  where $\abs{g^H}=\det((g^H)_{\lambda\mu})$.
It is easy to realize that $(j^1L)^*\cA$ is the Riemannian area
element on $L$ on every chart. In general, the above functional is
not global.
\begin{theorem}\label{th:arc-length}
  There are no global nowhere-vanishing Lagrangians on $J^1(E,1)$.
\end{theorem}
\begin{proof}
  If such a Lagrangian exists, then the line bundle $H^{1,0}$ would be
  orientable. But $H^{1,0}$ is never orientable. In fact, if $H^{1,0}$ would be
  orientable then its pull-back on any fiber of $J^1(E,1)$ would be orientable
  too. But this is a contradiction, because such a pull-back bundle is the
  canonical $1$-vector bundle over the Grassmannian manifold $\Gr (T_pE,1)$,
  that is not orientable~\cite{Osb82}.
\end{proof}
Note that we could get a global area Lagrangian if we would use
jets of \emph{oriented} submanifolds. Anyway we do not need to
postulate its globality in order to obtain global Euler--Lagrange
equations. Let us denote by $\Hess(\cA)$ the (local) Hessian of
$\cA$, \emph{i.e.}, the second differential of $\cA$ along the
fibers of $\pi_{1,0}$. We have
\begin{displaymath}
  \Hess(\cA)\colon J^1(E,n) \to V^*J^1(E,n)\otimes V^*J^1(E,n)\otimes
  \bigwedge^n (H^{1,0})^*,\quad
  \Hess(\cA)^{\lambda\mu}_{ij}=\pd{^2\sqrt{\abs{g^H}}}{u^i_\lambda\partial u^j_\mu}
\end{displaymath}
\begin{theorem}
  We have the equality
 \begin{displaymath}
    \cE(\cA) = - \mathbf{II}\circ \Hess(\cA),\quad\text{where}
    \quad \Hess(\cA)^{\lambda\mu}_{ij} =
    (\bar{g}^H)^{\lambda\mu}\,g^V_{ij}\sqrt{\abs{g^H}}.
 \end{displaymath}
\end{theorem}
\begin{proof}
Lemma \ref{le:vert_iso} and a long computation yields the result.
We used~\eqref{eq:E-L} and the formulas
  \begin{gather*}
    \pd{\sqrt{\abs{g^H}}}{u^i_\lambda} = \frac{1}{2} \sqrt{\abs{g^H}}
    (\bar{g}^H)^{\sigma\rho}\pd{(g^H)_{\sigma\rho}}{u^i_\lambda}, \quad
    \pd{(\bar{g}^H)^{\rho\theta}}{u^j_\mu} =
    -(\bar{g}^H)^{\rho\sigma}\pd{(g^H)_{\sigma\nu}}{u^j_\mu}(\bar{g}^H)^{\nu\theta}.
    \quad
  \end{gather*}
\end{proof}
\begin{corollary}
  The Euler--Lagrange equation $\cE(\cA)=0$ is an open submanifold
  of the minimal submanifold equation $\mathbf{H}^{-1}(0)$.
\end{corollary}

\subsection{Relativistic mechanics}\label{subsec:relat-mech}

Here we recall the geometric model for the phase space of one relativistic
particle by Jany\v ska and Modugno~\cite{JaMo96,Vit00}. The
phase space in this model is an open submanifold of the first-order jet of
curves in spacetime.  We will show that the equation of particle motion
presented in~\cite{JaMo96} is the Euler--La\-gran\-ge equation of
a Lagrangian, following the scheme of section~\ref{sec:var_princ}
(see~\cite{MaVi04} for more details).

\medskip

Let us set $\dim E=4$, and consider a \emph{scaled} Lorentz metric
\[
g\colon E\to \BL^2\otimes T^*E\otimes T^*E
\]
on $E$, with signature $(+---)$. Here $\BL^2$ is the one-dimensional space of
length units; we will also consider the mass and time one-dimensional spaces
$\BM$ and $\BT$ (see \cite{JaMoVi07} for more details on such spaces). The
components $g_{\mu\nu}$ of $g$ in local coordinates are $\BL^2$-valued smooth
functions on $E$. For physical reasons, we assume $E$ to be oriented and
time-like oriented. Latin indexes will label space-like coordinates, Greek
indexes will label spacetime coordinates. We will use coordinate charts
$(x^0,x^i)$ such that $\pd{}{x^0}$ is time-like and time-like oriented, and
$\pd{}{x^i}$ are space-like.

A \emph{motion} is a time-like curve $s\subset E$, and its \emph{velocity} is
$j_1s$. We consider the motion of a particle with mass $m\in\BM$. We introduce
the speed of light $c\in\BT^{-1}\otimes \BL$ and the Planck's constant
$\hbar\in\BT^{-1}\otimes\BL^2\otimes\BM$. The metric and the orientation yields
a natural parallelization $Ts\simeq s\times \BT$.

The set
\[
U^{1}E\subset J^1(E,1)
\]
of velocities of motions is said to be the {\em phase space}.  By a restriction
we have the natural projection $\pi_{1,0} \colon U^{1}E \to E$. Time
orientability implies that $H^{1,0}\simeq U^1E\times \BT$.

We introduce the normalized first-order contact map $\cyrm{D}\byd c D/\|D\|$,
with coordinate expression $\cyrm{D} = c\alpha\,D_0$, where $\alpha\byd (g_{00}
+ 2g_{0j}x^j_0 + g_{ij}x^i_0x^j_0)^{-1/2}$. The standard metric isomorphism
yields a natural $1$-form $\tau^\natural \byd g^\flat \circ \cyrm{D}$; in
coordinates $\tau^\natural = \tau_\lambda dx^\lambda = c\alpha (g_{0\lambda} +
g_{i\lambda}x^i_0) dx^\lambda$. Being $g (\cyrm{D},\cyrm{D}) = c^2$, the
inclusion $U^{1}E \subset \BT^*\otimes TE $ yields a non-linear bundle
structure on $U^1E$ with fibers diffeomorphic to $\R^3$.

After some intrinsic computations we obtain the \emph{gravitational $2$-form}
$\Omega^\natural$ on $U^1E$, with coordinate expression
\[
\Omega^\natural = c\alpha (g_{i\mu} - c^{-2} \tau_{i}\tau_{\mu})(dx^i_{0} -
\Gamma{_{\varphi}}{^i_{0}}) \wedge dx^{\mu},
\]
where $\Gamma_{\varphi}{^i_0} = K{_\varphi}{^i}_j x^j_0 + K{_\varphi}{^i}_0 -
x^i_0(K{_\varphi}{^0}_j x^j_0 + K{_\varphi}{^0}_0)$.  It can be proved that
$\Omega^\natural = d\tau^\natural$.

The \emph{electromagnetic field} can be introduced as a scaled closed form $F$
on $E$. It can be proved that the following \emph{joined contact $2$-form}
\[
\Omega \byd m/\hbar\,\Omega^\natural + q/(2\hbar c)\, F
\]
is non-degenerate.  Its kernel is a foliation $\gamma$ which yields the
dynamics on the spacetime.  If $A$ is a local potential of $F$ (according to
$2dA=F$), then $m/\hbar\,\tau^\natural+q/(\hbar c)\, A$ is a local potential of
$\Omega$.  Note that here $\hbar$ plays just the role of an overall scaling
factor, and it has no influence on the equation of motion.

We obtain our results following the scheme of section~\ref{sec:var_princ}, as
we did in the previous example.  First of all, we observe that
\[
\lambda_{GR}\byd \left[\frac{m}{\hbar}\tau^\natural+\frac{q}{\hbar c}A\right]=
h^{0,1}\left(\frac{m}{\hbar}\tau^\natural+\frac{q}{\hbar c}A\right)
\]
is a first-order local Lagrangian on $U^1E$, whose coordinate expression is
\[
\lambda_{GR}=\left(\frac{mc}{\hbar}\sqrt{g_{00} + 2g_{0j}x^j_0 +
    g_{ij}x^i_0x^j_0} + \frac{q}{\hbar c}\,(A_0+x^i_0A_i)\right)\,\odx^0.
\]
Then, the corresponding Euler--Lagrange morphism is
\[
\tilde{e}_1(\lambda_{GR})
=\left[d\left(\frac{m}{\hbar}\tau^\natural+\frac{q}{\hbar c}A\right)\right]
=[\Omega],
\]
with coordinate expression
\begin{equation}\label{eq:2009}
  \tilde{e}_1(\lambda_{GR})=\frac{mc}{\hbar}\alpha (g_{ij} - c^{-2}
  \tau_i \tau_j)(x^i_{00}-
  (\gamma^i_{00}{}^\natural+\gamma^i_{00}{}^e))\, \omega^j\otimes
  \odx^0,
\end{equation}
where $\gamma^i_{00}{}^\natural+\gamma^i_{00}{}^e$ are the components of the
above foliation $\gamma$, with coordinate expression
\begin{align*}
  & \gamma^i_{00}{}^\natural=K{_0}{^i}{_0}-2K{_0}{^i}{_j}x^j_0+
  K{_0}{^0}{_0}x^i_0+2K{_0}{^0}{_j}x^i_0x^j_0-K{_j}{^i}{_k}x^j_0x^k_0
  +K{_j}{^0}{_k}x^j_0x^k_0x^i_0,
  \\
  &\gamma^i_{00}{}^e= - \frac{q}{mc}(g^{i\mu} - x^i_0 g^{0\mu})(F_{0\mu} +
  F_{j\mu}x^j_0).
\end{align*}
Eq.~(\ref{eq:2009}) is equivalent to
$x^i_{00}-(\gamma^i_{00}{}^\natural+\gamma^i_{00}{}^e)=0$, and coincide with
the equation of the integral curves of the foliation $\gamma$. This is the
equation of particle motion in general relativity, written for non-parametrized
time-like curves.  Note that $\lambda_{GR}$ has the same domain as $A$, hence
it is global if and only if $A$ is global (or vanishing).
\begin{remark}
  The relativistic mechanics on jets of submanifolds has a distinguished
  feature: it can be easily proved that its Lagrangian is non-degenerate. In
  other words, $U^1E$ is a natural first class Dirac constraint (see
  \cite{MaVi04} for more details).
\end{remark}

%-----------------------------------------------------------------------------%
% Section
%-----------------------------------------------------------------------------%
\section{Conclusions}
\label{sec:conclusions}
The fact that the objects in the $\cC$-spectral sequence on submanifolds can be
uniquely represented through the bundle of total derivatives allows us to
consider a number of problems by analogy with the fibered case. Among such
problems we have the classification of Lagrangians and/or Euler--Lagrange forms
of a given order which are invariant under the action of a symmetry group or
pseudogroup. In the Riemannian case, a first approach to the above problem is
in \cite{MuMu04}, where isometry-invariant Lagrangians on jets of
  immersions are classified. Our language allows one to formulate the
problem for more general objects, like Euler--Lagrange morphisms.

\medskip
\textbf{Acknowledgements:}

\smallskip
  We thank D. Canarutto, M. Modugno, D. Krupka, O. Krupkov\'a, A.M.
  Verbovetsky for many stimulating discussions that helped us to improve the
  exposition and the results. The first author would also like to thank S.K.
  Donaldson and D.J. Saunders for having read and commented parts of the
  material in this paper for the discussion of his Ph.D. thesis
  \cite{MaThesis}.

%----------------------------------------------------------------------------%
% B I B L I O G R A P H Y
%----------------------------------------------------------------------------%

\end{document}